\documentclass{amsart}

\newenvironment{thmenum}{

\begin{enumerate}}
{\end{enumerate}}

\usepackage[dvips]{graphicx}
\usepackage{epsfig, latexsym, amsfonts, amsthm, amsmath, pstricks, mathtools}

\renewenvironment{proof}[1][]{\vskip-\lastskip\par\vskip6pt plus2pt
minus0pt\par%
\noindent\textit{Proof.}\enspace\ignorespaces}{\hfill$\Box$\par\vskip6pt
plus2pt minus0pt}

\numberwithin{equation}{section}

\newtheorem{theorem}[equation]{Theorem}

\newtheorem{lemma}[equation]{Lemma}

\newtheorem{proposition}[equation]{Proposition}

\newtheorem{corollary}[equation]{Corollary}

\theoremstyle{definition}

\newtheorem{notation}[equation]{Notation}

\begin{document}
\title{The norm of a Ree Group}
\author{Tom De Medts}
\address{Department of Pure Mathematics and Computer Algebra \\
	Ghent University \\
	9000 Gent, Belgium}
\email{tdemedts@cage.ugent.be}
\author{Richard M. Weiss}
\address{Department of Mathematics \\
         Tufts University \\
         Medford, MA 02155, USA}
\email{rweiss@tufts.edu}


\date{June 2, 2009}

\begin{abstract}
We give an explicit construction of the Ree groups of type $G_2$ as 
groups acting on mixed Moufang hexagons together with detailed proofs
of the basic properties of these groups contained in the two fundamental
papers of Tits on this subject, \cite{tits} and \cite{octagons}. We also give a 
short proof that the norm of a Ree group is anisotropic.
\end{abstract}

\maketitle

\section{Introduction}

The finite Ree groups of type $G_2$ were introduced by Ree in \cite{ree}.
In \cite{tits}, Tits showed how to construct these groups over an arbitrary
field $K$ in characteristic three having an endomorphism 
whose square is the Frobenius endomorphism of $K$ and gave a number
of their basic properties. 
His result can be summarized as follows.

\begin{theorem}\label{ree1}
Let $K$ be a field of characteristic 3 and suppose that $K$ has
an endomorphism $\theta$ such that 
\[ x^{\theta^2}=x^3 \]
for all $x\in K$. Let $U$ denote the set $K\times K\times K$ endowed with the 
multiplication:
\begin{equation}\label{ree1a}
(a,b,c)\cdot(x,y,z)=(a+x,b+y+a^\theta x,c+z-ay+bx-a^{\theta+1}x) \,,
\end{equation}
and let 
\begin{equation}\label{ree1b}
H=\{h_t\mid t\in K^*\} \,,
\end{equation}
where for each $t\in K^*$, $h_t$ is the map from $U$ to itself given by the formula
\[ (a,b,c)^{h_t}=(ta,t^{\theta+1}b,t^{\theta+2}c) \,. \]
Let
\begin{equation}\label{ree1c}
N(a,b,c)=-ac^\theta+a^{\theta+1}b^\theta-a^{\theta+3}b-a^2b^2+b^{\theta+1}+c^2-a^{2\theta+4}
\end{equation}
for all $(a,b,c)\in U$ and let $X$ denote the disjoint union of $U$ and a symbol $\infty$.
Then the following hold:
\begin{thmenum}
\item $U$ is a group with identity $(0,0,0)$ (which we denote by $0$) and inverses given by
\[ (a,b,c)^{-1}=(-a,-b+a^{\theta+1},-c) \]
and $H$ is a
group of automorphisms of $U$. 
\item The map $N$ is anisotropic. This is to say, $N(a,b,c)=0$ if and only if $(a,b,c)=0$.
\item Let $\omega$ be the map from $X$ to itself that interchanges $\infty$ and $0$ and maps
an arbitrary element $(a,b,c)$ of $U^*$ to 
\begin{equation}\label{ree1d}
\big(-v/w,-u/w,-c/w),
\end{equation}
where $v=a^\theta b^\theta-c^\theta+ab^2+bc-a^{2\theta+3}$, 
$u=a^2b-ac+b^\theta-a^{\theta+3}$ and $w=N(a,b,c)$.
Let $U$ be identified with the permutation group of $X$ that fixes $\infty$ and acts
on $X\backslash\{\infty\}$ by right multiplication. Let $H$ be identified
with the permutation group of $X$ that fixes $\infty$ and acts on $X\backslash\{\infty\}$
by the formula~\eqref{ree1a} (and thus fixes also $0$). 
Let $K^\dagger$ be the subgroup of $K^*$ generated by $\{N(a,b,c)\mid (a,b,c)\in U^*\}$
and let 
\begin{equation}\label{ree1e}
H^\dagger=\{h_t\mid t\in K^\dagger\}\subset H.
\end{equation}
Then $\omega$ is a permutation of $X$ of order~2 and the
subgroup $G$ of ${\rm Sym}(X)$ generated by $U$ and $\omega$ has the following properties:
\begin{thmenum}
\item[(I)] $G$ is a 2-transitive permutation group on $X$.
\item[(II)] $U$ is a normal subgroup of the stabilizer $G_\infty$ and $G_\infty =UH^\dagger$.
\item[(III)] $G=\langle U,U^\omega\rangle$.
\item[(IV)] $H$ normalizes $G$.
\item[(V)] $\omega$ inverts every element of $H$.
\item[(VI)] If $|K|>3$, then $G$ is simple.
\end{thmenum}
\end{thmenum}
\end{theorem}

Tits' proof of Theorem~\ref{ree1} in \cite{tits} is based on the 
standard embedding of the split Moufang hexagon in 6-dimensional
projective space; see also \cite[Section~7.7]{hendrik}. 
The purpose of this note is to give an alternative proof of Theorem~\ref{ree1}
in which we construct the set $X$ inside the mixed hexagon
defined over the pair $(K,K^\theta)$, which we construct directly
without reference to projective space.

Our motivation is threefold. First, since the Ree groups of type $G_2$
continue to be the center of lively interest (see especially
\cite{guralnick}), we want give a proof of Theorem~\ref{ree1}
in which many of the details left to the
reader in \cite{tits} are filled in. We also want to provide independent confirmation
of the accuracy of the formulas occurring in Theorem~\ref{ree1}.
(In fact, there is a $\theta$ missing 
in the second term in the definition of the norm and a minus sign missing 
in front of the whole expression on page~12 in \cite{tits}, where
$\theta$ is called $\sigma$ and the norm $N$ is called $w$.)
Secondly, we want to examine 
the fact that the map $N$, which we
call the {\it norm} of $G$, is anisotropic. As in \cite{tits},
this fact emerges ``geometrically'' in the course of our proof
of Theorem~\ref{ree1}; in Section~\ref{ree11}, we give
a succinct algebraic explanation. Thirdly, we hope that the method we use to
prove Theorem~\ref{ree1} 
can serve as a model for other calculations in Moufang
polygons and in more general types of buildings.

If $|K|=3$, then the endomorphism $\theta$ is trivial and
the group $G$ is not simple; in fact, it is isomorphic to 
${\rm Aut}(L_2(8))$ in this case and thus has a normal subgroup of index~3 (which is
simple).

If $K$ is finite, then $H^\dagger=H$ and thus $H\subset G$ 
(by \cite[8.4]{ree}). It is not true
in general, however, that $H=H^\dagger$.
We say a few words about this in Section~\ref{ree18}.

We mention that there are also Ree groups of type $F_4$. The canonical
reference for these groups is \cite{octagons}.

\smallskip
\section{The hexagon of mixed type}

Let $K$ be a field of characteristic~3 and let
$\theta$ be a square root of the Frobenius endomorphism 
of $K$. We now begin our proof of Theorem~\ref{ree1} by 
constructing the mixed hexagon associated with the pair $(K,\theta)$.
(See \cite[16.20 and 41.20]{TW}
for the definition of a mixed hexagon.)
Let $U_1,U_2,\ldots,U_6$
be six groups isomorphic to the additive group of $K$
and let for each $i\in [1,6]$, let $x_i$ be an isomorphism from 
$K$ to $U_i$. Let $U_+$ be the group generated
by the groups $U_1,U_2,\ldots,U_6$ subject to the commutator relations
\begin{align}
[x_1(s),x_5(t)]&=x_3(-st)\notag\\
[x_2(s),x_6(t)]&=x_4(st)\text{ and}\label{ree10}\\
[x_1(s),x_6(t)]=x_2(-s^\theta t)&x_3(-s^2 t^\theta)x_4(s^\theta t^2)x_5(st^\theta)\notag
\end{align}
for all $s,t\in K$ and $[U_i,U_j]=1$ for all other pairs $i,j$ such that $1\le i<j\le 6$. 
(We are using the convention that $[a,b]=a^{-1}b^{-1}ab=(b^{-1})^ab$.)
Every element of $U_+$ can be written uniquely as an element in 
the product $U_1U_2\cdots U_6$ and there is an automorphism $\rho$
of $U_+$ interchanging $x_i(t)$ and $x_{7-i}(t)$ for all
$i\in[1,6]$ and all $t\in K$. We will see below that the group $U$ 
in Theorem~\ref{ree1} is the centralizer of $\rho$ in $U_+$.

Let $U_{i,j}$ denote the subgroup $U_iU_{i+1}\cdots U_j$ of $U_+$ 
for all $i,j$ such that $1\le i\le j\le 6$ (so $U_{i,i}=U_i$ for each $i$).
For each $i\in[1,5]$, let $W_i$ denote the set of right cosets in $U_+$ of $U_{1,6-i}$.
For each $i\in[6,10]$, let $W_i$ denote the set of right cosets in $U_+$ of
$U_{12-i,6}$. 
Let $W$ be the disjoint union of $W_1,W_2,\ldots,W_{10}$ together with two symbols
$\bullet$ and $\star$. For each $i\in[1,9]$, let $E_i$ be the set of pairs
$\{x,y\}$ such that $x\in W_i$, $y\in W_{i+1}$ and the intersection of $x$ and $y$
is non-empty. Let $E$ be the set of (unordered) 2-element subsets of $W$
consisting of $\{\bullet,\star\}$, $\{\bullet,x\}$ for all $x\in W_1$,
$\{\star,y\}$ for all $y\in W_{10}$ together with all the pairs in 
$E_1\cup E_2\cup\ldots\cup E_9$. Finally, let $\Gamma$ be the graph with vertex set 
$W$ and edge set $E$. 

\setlength{\unitlength}{1cm}
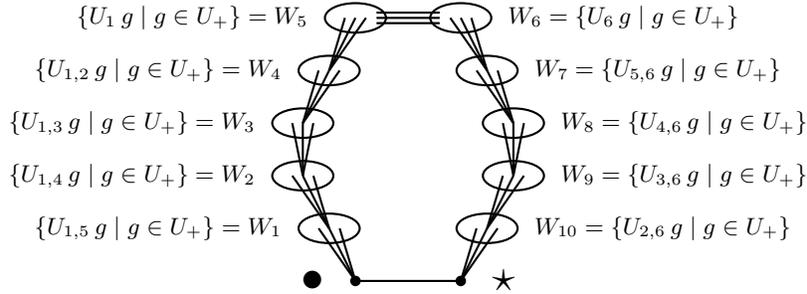
\begin{figure}[ht!]
\begin{center}
\psset{unit=.7}
\begin{pspicture}(-5,0)(7,5.3)
	\pscircle*(0,0){.1}
	\pscircle*(2,0){.1}
	\psellipse(2.5,1)(.6,.3)
	\psellipse(3,2)(.6,.3)
	\psellipse(3,3)(.6,.3)
	\psellipse(2.5,4)(.6,.3)
	\psellipse(2,5)(.6,.3)
	\psellipse(-.5,1)(.6,.3)
	\psellipse(-1,2)(.6,.3)
	\psellipse(-1,3)(.6,.3)
	\psellipse(-.5,4)(.6,.3)
	\psellipse(0,5)(.6,.3)
	\psline(0,0)(2,0)
	\psline(0,0)(-.7,1) \psline(0,0)(-.5,1) \psline(0,0)(-.3,1)
	\psline(-.5,1)(-1.2,2) \psline(-.5,1)(-1.0,2) \psline(-.5,1)(-0.8,2)
	\psline(-1,2)(-1.2,3) \psline(-1,2)(-1.0,3) \psline(-1,2)(-0.8,3)
	\psline(-1,3)(-.7,4) \psline(-1,3)(-.5,4) \psline(-1,3)(-.3,4)
	\psline(-.5,4)(-.2,5) \psline(-.5,4)(0,5) \psline(-.5,4)(.2,5)
	\psline(2,0)(2.3,1) \psline(2,0)(2.5,1) \psline(2,0)(2.7,1)
	\psline(2.5,1)(2.8,2) \psline(2.5,1)(3.0,2) \psline(2.5,1)(3.2,2)
	\psline(3,2)(2.8,3) \psline(3,2)(3.0,3) \psline(3,2)(3.2,3)
	\psline(3,3)(2.3,4) \psline(3,3)(2.5,4) \psline(3,3)(2.7,4)
	\psline(2.5,4)(1.8,5) \psline(2.5,4)(2,5) \psline(2.5,4)(2.2,5)
	\psline(0.4,5)(1.6,5)
	\psline(0.4,5.1)(1.6,5.1)
	\psline(0.4,4.9)(1.6,4.9)
	\rput[r](-.6,0){\huge $\bullet$}
	\rput[l](2.6,0){\huge $\star$}
	\rput[r](-1.4,1){\small $\{ U_{1,5}\,g \mid g \in U_+ \} = W_1$}
	\rput[l](3.4,1){\small $W_{10} = \{ U_{2,6}\,g \mid g \in U_+ \}$}
	\rput[r](-1.9,2){\small $\{ U_{1,4}\,g \mid g \in U_+ \} = W_2$}
	\rput[l](3.9,2){\small $W_9 = \{ U_{3,6}\,g \mid g \in U_+ \}$}
	\rput[r](-1.9,3){\small $\{ U_{1,3}\,g \mid g \in U_+ \} = W_3$}
	\rput[l](3.9,3){\small $W_8 = \{ U_{4,6}\,g \mid g \in U_+ \}$}
	\rput[r](-1.4,4){\small $\{ U_{1,2}\,g \mid g \in U_+ \} = W_4$}
	\rput[l](3.4,4){\small $W_7 = \{ U_{5,6}\,g \mid g \in U_+ \}$}
	\rput[r](-0.9,5){\small $\{ U_{1}\,g \mid g \in U_+ \} = W_5$}
	\rput[l](2.9,5){\small $W_6 = \{ U_{6}\,g \mid g \in U_+ \}$}
\end{pspicture}
\caption{The graph $\Gamma$ \label{fig1} }
\end{center}
\end{figure}

\begin{proposition}\label{ree2}
The graph $\Gamma$ is the Moufang hexagon associated with the hexagonal system
$(K/K^\theta)^\circ$ as defined in \cite[15.20 and 16.8]{TW}.
\end{proposition}

\begin{proof}
The maps $x_i(s)\mapsto x_i(s^\theta)$ for $i=2$, $4$ and $6$,
$x_i(s)\mapsto x_i(-s)$ for $i=3$ and $5$
and $x_1(s)\mapsto x_1(s)$
extend to an isomorphism from the group $U_+$ defined
above to the group 
obtained by setting $F=K^\theta$, $J=K$, $T(a,b)=0$, $a^\#=a^2$,
$N(a)=a^3$ and $a\times b=2ab$ for all $a,b\in K$ in 
\cite[16.8]{TW}. By \cite[7.5]{TW},
it follows that the graph $\Gamma$ is precisely the
Moufang hexagon associated with the hexagonal system
$(K/K^\theta)^\circ$.
\end{proof}

\begin{notation}\label{ree16}
Let $D={\rm Aut}(\Gamma)$ and let $D^\dagger$ denote the subgroup of 
$D$ generated by all the root groups of $\Gamma$.
\end{notation}

The group $U_+$ acts faithfully by right
multiplication on the elements of
\[ W_1\cup\cdots\cup W_{10} \]
and maps the set $E$ of edges of $\Gamma$ to itself.
This allows us to identify $U_+$ with a subgroup of $D^\dagger$
fixing the two vertices $\bullet$ and $\star$. Thus, for example,
we have
\begin{equation}\label{ree16a}
U_{15}^{x_6(t)}=U_{15}x_6(t),
\end{equation}
where the cosets $U_{15}$ and $U_{15}x_6(t)$ are vertices in the set $W_1$
and the expression on the left means the image of the vertex $U_{15}$ under
the action of the element $x_6(t)\in U_+$.

We observe that by \cite[35.13]{TW}, $\Gamma$ is not a split hexagon
unless $K=K^\theta$. In other words, $\Gamma$ is split if and only
if $K$ is perfect.

\smallskip
\section{The automorphisms $m_1$ and $m_6$}\label{ree17}

Let $\Sigma$ denote the apartment of $\Gamma$
spanned by the vertices $\bullet$, $\star$, $U_{1,6-i}\in W_i$ for
all $i\in[1,5]$ and $U_{12-i,6}\in W_i$ for all $i\in[6,10]$. 
From now on, we will write $U_{ij}$ in place of $U_{i,j}$.
Let 
\[ m_1=\mu(x_1(1))\text{ and }m_6=\mu(x_6(1)) \,, \]
where the map $\mu$ is defined (with respect to the apartment $\Sigma$)
as in \cite[6.1]{TW}. 
Both of these elements are contained in the group $D^\dagger$ and 
both induce reflections on $\Sigma$, 
$m_1$ the reflection fixing $\star$ and $U_1$ and $m_6$ the reflection
fixing $\bullet$ and $U_6$.
By \cite[32.12]{TW}, we have
\[ x_6(t)^{m_1}=x_2(t)\text{ and }x_5(t)^{m_1}=x_3(t) \]
and
\[ x_1(t)^{m_6}=x_5(-t)\text{ and }x_2(t)^{m_6}=x_4(t) \]
for all $t\in K$. Thus the action of $m_1$ on the vertices in $W_1$ is given by
\begin{equation}\label{ree4a}
(U_{15}x_6(t))^{m_1}=U_{15}^{x_6(t)m_1}=U_{15}^{m_1x_2(t)}=U_{36}x_2(t)
\end{equation}
-- see (\ref{ree16a}) above -- and the action of $m_6$ on the vertices in $W_{10}$ is given by
\begin{equation}\label{ree4b}
(U_{26}x_1(t))^{m_6}=U_{26}^{x_1(t)m_6}=U_{26}^{m_6x_5(-t)}=U_{14}x_5(-t)
\end{equation}
for all $t\in K$. Similarly, we have
\begin{equation}\label{ree4c}
(U_{14}x_5(t))^{m_1}=U_{46}x_3(t)
\end{equation}
and 
\begin{equation}\label{ree4d}
(U_{36}x_2(t))^{m_6}=U_{13}x_4(t)
\end{equation}
for all $t\in K$.

\begin{proposition}\label{ree3}
The maps $m_1$ and $m_6$ are as in Tables \ref{ta:m1} and \ref{ta:m6}.
(For use in Section~\ref{ree15}, we have recorded also the product
$m_1m_6$ in Table \ref{ta:m16}.)
\end{proposition}

\begin{proof}
It is a simple calculation to check using the commutator relations \eqref{ree10} that the permutation of
the vertex set $W$ described in Table \ref{ta:m1} is an automorphism
of the hexagon $\Gamma$. Moreover, this automorphism
induces the same reflection on $\Gamma$ as does $m_1$ and 
it agrees with $m_1$ on the set of neighbors of $\bullet$
and on the set of neighbors of $U_{15}$ by \eqref{ree4a} and \eqref{ree4c}.
By \cite[3.7]{TW}, it follows that this automorphism equals $m_1$.
By \eqref{ree4b}, \eqref{ree4d} and a similar argument,
the claim holds for $m_6$.
\end{proof}

\begin{table}
\small
\begin{align*}
	\star &\mapsto \star \\
	\bullet &\mapsto U_{26} \\
	U_{15}\,x_6(t) &\mapsto U_{36}\,x_2(t) \\
	U_{14}\,x_5(s)\,x_6(t) &\mapsto U_{46}\,x_2(t)\,x_3(s) \\
	U_{13}\,x_4(r)\,x_5(s)\,x_6(t) &\mapsto U_{56}\,x_2(t)\,x_3(s)\,x_4(r) \\
	U_{12}\,x_3(u)\,x_4(r)\,x_5(s)\,x_6(t) &\mapsto U_{6}\,x_2(t)\,x_3(s)\,x_4(r)\,x_5(-u) \\
	U_{1}\,x_2(v)\,x_3(u)\,x_4(r)\,x_5(s)\,x_6(t) &\mapsto U_{1}\,x_2(t)\,x_3(s)\,x_4(r+vt)\,x_5(-u)\,x_6(-v) \\[1ex]
	U_{6}\,x_1(s)\,x_2(t)\,x_3(r)\,x_4(u)\,x_5(v) &\xmapsto{s = 0} U_{12}\,x_3(v)\,x_4(u)\,x_5(-r)\,x_6(-t) \\[-1.2ex]
		&\xmapsto{s \neq 0} U_{6}\,x_1(-s^{-1})\,x_2(-s^{-\theta}t)\,x_3(v + s^{-2}t^\theta) \\[-.6ex]
		&\hspace*{12ex} \cdot x_4(u - s^{-\theta}t^2)\,x_5(s^{-1}t^\theta - r) \\
	U_{56}\,x_1(s)\,x_2(t)\,x_3(r)\,x_4(u) &\xmapsto{s = 0} U_{13}\,x_4(u)\,x_5(-r)\,x_6(-t) \\[-1.2ex]
		&\xmapsto{s \neq 0} U_{56}\,x_1(-s^{-1})\,x_2(-s^{-\theta}t)\,x_3(-s^{-1}r - s^{-2}t^\theta) \\[-.6ex]
		&\hspace*{12ex} \cdot x_4(u - s^{-\theta}t^2) \\
	U_{46}\,x_1(s)\,x_2(t)\,x_3(r) &\xmapsto{s = 0} U_{14}\,x_5(-r)\,x_6(-t) \\[-1.2ex]
		&\xmapsto{s \neq 0} U_{46}\,x_1(-s^{-1})\,x_2(-s^{-\theta}t)\,x_3(-s^{-1}r - s^{-2}t^\theta) \\
	U_{36}\,x_1(s)\,x_2(t) &\xmapsto{s = 0} U_{15}\,x_6(-t) \\[-1.2ex]
		&\xmapsto{s \neq 0} U_{36}\,x_1(-s^{-1})\,x_2(-s^{-\theta}t) \\
	U_{26}\,x_1(s) &\xmapsto{s = 0} \bullet \\[-1.2ex]
		&\xmapsto{s \neq 0} U_{26}\,x_1(-s^{-1})
\end{align*}
\caption{The action of $m_1$ on $\Gamma$\label{ta:m1}}
\end{table}

\begin{table}
\small
\begin{align*}
	\star &\mapsto U_{15} \\
	\bullet &\mapsto \bullet \\
	U_{15}\,x_6(w) &\xmapsto{w = 0} \star \\[-1.2ex]
		&\xmapsto{w \neq 0} U_{15}\,x_6(-w^{-1}) \\
	U_{14}\,x_5(v)\,x_6(w) &\xmapsto{w = 0} U_{26}\,x_1(-v) \\[-1.2ex]
		&\xmapsto{w \neq 0} U_{14}\,x_5(-vw^{-\theta})\,x_6(-w^{-1}) \\
	U_{13}\,x_4(u)\,x_5(v)\,x_6(w) &\xmapsto{w = 0} U_{36}\,x_1(-v)\,x_2(-u) \\[-1.2ex]
		&\xmapsto{w \neq 0} U_{13}\,x_4(-v^\theta w^{-2} - w^{-1}u)\,x_5(-vw^{-\theta})\,x_6(-w^{-1}) \\
	U_{12}\,x_3(r)\,x_4(u)\,x_5(v)\,x_6(w) &\xmapsto{w = 0} U_{46}\,x_1(-v)\,x_2(-u)\,x_3(r) \\[-1.2ex]
		&\xmapsto{w \neq 0} U_{12}\,x_3(r - v^2 w^{-\theta})\,x_4(-v^\theta w^{-2} - w^{-1}u) \\[-.6ex]
		&\hspace*{10ex} \cdot x_5(-vw^{-\theta})\,x_6(-w^{-1}) \\
	U_{1}\,x_2(t)\,x_3(r)\,x_4(u)\,x_5(v)\,x_6(w) &\xmapsto{w = 0} U_{56}\,x_1(-v)\,x_2(-u)\,x_3(r)\,x_4(t) \\[-1.2ex]
		&\xmapsto{w \neq 0} U_{1}\,x_2(v^\theta w^{-1} - u - tw)\,x_3(r - v^2 w^{-\theta}) \\[-.6ex]
		&\hspace*{10ex} \cdot x_4(-v^\theta w^{-2} - w^{-1}u)\,x_5(-vw^{-\theta})\,x_6(-w^{-1}) \\
	U_{6}\,x_1(s)\,x_2(t)\,x_3(r)\,x_4(u)\,x_5(v) &\mapsto U_{6}\,x_1(-v)\,x_2(-u)\,x_3(r-sv)\,x_4(t)\,x_5(s) \\[1ex]
	U_{56}\,x_1(s)\,x_2(t)\,x_3(r)\,x_4(u) &\mapsto U_{1}\,x_2(-u)\,x_3(r)\,x_4(t)\,x_5(s) \\
	U_{46}\,x_1(s)\,x_2(t)\,x_3(r) &\mapsto U_{12}\,x_3(r)\,x_4(t)\,x_5(s) \\
	U_{36}\,x_1(s)\,x_2(t) &\mapsto U_{13}\,x_4(t)\,x_5(s) \\
	U_{26}\,x_1(s) &\mapsto U_{14}\,x_5(s)
\end{align*}
\caption{The action of $m_6$ on $\Gamma$\label{ta:m6}}
\end{table}

\begin{table}
\small
\begin{align*}
	\star &\mapsto U_{15} \\
	\bullet &\mapsto U_{14} \\
	U_{15}\,x_6(t) &\mapsto U_{13}\,x_4(t) \\
	U_{14}\,x_5(s)\,x_6(t) &\mapsto U_{12}\,x_3(s)\,x_4(t) \\
	U_{13}\,x_4(r)\,x_5(s)\,x_6(t) &\mapsto U_{1}\,x_2(-r)\,x_3(s)\,x_4(t) \\
	U_{12}\,x_3(u)\,x_4(r)\,x_5(s)\,x_6(t) &\mapsto U_{6}\,x_1(u)\,x_2(-r)\,x_3(s)\,x_4(t) \\[1ex]
	U_{1}\,x_2(v)\,x_3(u)\,x_4(r)\,x_5(s)\,x_6(t) &\xmapsto{v = 0} U_{56}\,x_1(u)\,x_2(-r)\,x_3(s)\,x_4(t) \\[-1.2ex]
		&\xmapsto{v \neq 0} U_{1}\,x_2(u^\theta v^{-1} - r) \, x_3(s + u^2 v^{-\theta}) \\[-.6ex]
		&\hspace*{12ex} \cdot x_4(u^\theta v^{-2} + v^{-1}r + t) \, x_5(-u v^{-\theta}) \, x_6(v^{-1}) \\
	U_{6}\,x_1(s)\,x_2(t)\,x_3(r)\,x_4(u)\,x_5(v) &\xmapsto{s = 0, t = 0} U_{46}\,x_1(r)\,x_2(-u)\,x_3(v) \\[-1.2ex]
		&\xmapsto{s = 0, t \neq 0} U_{12}\,x_3(v + r^2 t^{-\theta})\,x_4(r^\theta t^{-2} + t^{-1}u) \\[-.6ex]
		&\hspace*{12ex} \cdot x_5(-r t^{-\theta})\,x_6(t^{-1}) \\[-1.2ex]
		&\xmapsto{s \neq 0} U_{6}\,x_1(r - s^{-1} t^\theta)\,x_2(s^{-\theta}t^2 - u) \\[-.6ex]
		&\hspace*{12ex} \cdot x_3(v - s^{-2} t^\theta - s^{-1} r) \, x_4(-s^{-\theta} t) \, x_5(-s^{-1}) \\
	U_{56}\,x_1(s)\,x_2(t)\,x_3(r)\,x_4(u) &\xmapsto{s = 0, t = 0} U_{36}\,x_1(r)\,x_2(-u) \\[-1.2ex]
		&\xmapsto{s = 0, t \neq 0} U_{13}\,x_4(r^\theta t^{-2} + t^{-1}u) \, x_5(-rt^{-\theta}) \, x_6(t^{-1}) \\[-.6ex]
		&\xmapsto{s \neq 0} U_{1}\,x_2(-u + s^{-\theta}t^2)\,x_3(-s^{-1}r - s^{-2}t^\theta) \\[-.6ex]
		&\hspace*{12ex} \cdot x_4(- s^{-\theta}t) \, x_5(-s^{-1}) \\
	U_{46}\,x_1(s)\,x_2(t)\,x_3(r) &\xmapsto{s = 0, t = 0} U_{26}\,x_1(r) \\[-1.2ex]
		&\xmapsto{s = 0, t \neq 0} U_{14}\,x_5(-rt^{-\theta}) \, x_6(t^{-1}) \\[-1.2ex]
		&\xmapsto{s \neq 0} U_{12}\,x_3(-s^{-1}r - s^{-2}t^\theta) \, x_4(-s^{-\theta}t) \, x_5(-s^{-1}) \\
	U_{36}\,x_1(s)\,x_2(t) &\xmapsto{s = 0, t = 0} \star \\[-1.2ex]
		&\xmapsto{s = 0, t \neq 0} U_{15}\,x_6(t^{-1}) \\[-1.2ex]
		&\xmapsto{s \neq 0} U_{13}\,x_4(-s^{-\theta}t) \, x_5(-s^{-1}) \\
	U_{26}\,x_1(s) &\xmapsto{s = 0} \bullet \\[-1.2ex]
		&\xmapsto{s \neq 0} U_{14}\,x_5(-s^{-1})
\end{align*}
\caption{The action of $m_1 m_6$ on $\Gamma$\label{ta:m16}}
\end{table}

Now let $\rho$ be the automorphism of $U_+$ mentioned above.
Thus 
\begin{equation}\label{ree23}
x_i(t)^\rho=x_{7-i}(t)
\end{equation} 
for all $i\in[1,6]$ and all
$t\in K$. By \cite[7.5]{TW}, $\rho$ extends to 
a unique automorphism of $\Gamma$ which we also denote
by $\rho$. The automorphism $\rho$ induces the reflection
on the apartment $\Sigma$ that interchanges $\bullet$ and 
$\star$ as well as $U_1$ and $U_6$. 

From now on, we set 
\begin{equation}\label{ree9}
\omega=(m_1m_6)^3 \,.
\end{equation}

\begin{proposition}\label{ree5}
The automorphisms $\rho$ and $\omega$ commute with each
other and both have order~2.
\end{proposition}

\begin{proof}
Since $\rho$ has order~2 as an automorphism of $U_+$, it
also has order~2 as an automorphism of $\Gamma$. 
By \cite[6.9]{TW}, $\omega=(m_6m_1)^3$ and by
\cite[6.2]{TW}, $m_1^\rho=m_6$ and $m_6^\rho=m_1$. 
Thus $\omega^\rho=(m_6m_1)^3=\omega$. Let
$d=m_1^2$ and $e=m_6^2$ (so $[m_1,d]=[m_6,e]=1$). Then $d$ and $e$ both act
trivially on the apartment $\Sigma$ and by \cite[29.12]{TW}, 
$d$ centralizes $U_1$ and $U_4$ and inverts every 
element of $U_i$ for 
all other $i\in[1,6]$ and $e$ centralizes $U_3$
and $U_6$ and inverts every element of $U_i$ for all
other $i\in[1,6]$. By \cite[6.7]{TW}, $d$ and $e$ are elements of order~2
(so $m_1^{-1}=dm_1$ and $m_6^{-1}=em_6$)
and their product (in either order) is the unique element of $D$
acting trivially on $\Sigma$ that centralizes $U_2$ and $U_5$
and inverts every element of $U_i$ for all other $i\in[1,6]$.
Since
$U_i^{m_1}=U_{8-i}$ for all $i\in[2,6]$ and
$U_i^{m_6}=U_{6-i}$ for all $i\in[1,5]$, both 
$e^{m_1}$ and $d^{m_6}$ centralize $U_2$ and $U_5$ and 
invert every element of $U_i$ for all other $i\in[1,6]$.
Thus $e^{m_1}=ed=d^{m_6}$. It follows by repeated use of these relations that
\[ (m_1^{-1}m_6^{-1})^3=(dm_1\cdot em_6)^3=(m_1m_6)^3 \]
and hence $\omega^{-1}=(m_6m_1)^{-3}=\omega$.
\end{proof}

\begin{proposition}\label{ree6}
Let $\varphi$ be the map from $U$ to $U_+$ given by 
\[ \varphi(a,b,c) = x_1(a) \, x_2(b) \, x_3(c - ab + a^{\theta+2}) \, x_4(c + ab) \,
	x_5(b - a^{\theta+1}) \, x_6(a) \,. \]
Then $\varphi$ is an injective homomorphism whose image
is the centralizer of $\rho$ in $U_+$.
\end{proposition}

\begin{proof}
This holds by \eqref{ree1a} and \eqref{ree10}. 
\end{proof}

\noindent
From now on we identify $U$ with its image in $U_+$ under the map $\varphi$ in Proposition~\ref{ree6}.

\begin{proposition}\label{ree7}
Let $X$ be the set of edges of $\Gamma$ fixed by $\rho$, let $\infty$ denote
the edge $\{\bullet,\star\}$ and 
let $G=\langle U,\omega\rangle$, where $\omega$ is as in \eqref{ree9}.
Then the following hold:
\begin{thmenum}
\item $U$ acts regularly on $X\backslash\{\infty\}$.
\item $G$ acts 2-transitively on $X$. 
\item $G=B\cup B\omega B$, where $B=G_\infty$. 
\item $U$ is a normal subgroup of the stabilizer $G_\infty$.
\item $G$ acts faithfully on $X$. 
\end{thmenum}
\end{proposition}

\begin{proof}
Since $\rho$ interchanges the vertices $\bullet$ and $\star$, all the 
edges in $X$ other than $\infty=\{\bullet,\star\}$ are two-element
subsets containing a
right coset of $U_1$ and a right coset of $U_6$. Since every element
of $U_+$ has a unique representation in $U_1U_2\cdots U_6$, the
intersection of a right
coset of $U_1$ and a right coset of $U_6$ is either empty or consists
of a unique element. It follows that
\[ X=\big\{\{U_1g,U_6g\}\mid g\in U\big\}\cup\{\infty\} \,. \]
In particular, (i) holds 
and we can identify $U$ with $X\backslash\{\infty\}$ via
the map that sends $g\in U$ to $\{U_1,U_6\}^g$. In
particular, $0$ now denotes the edge $\{U_1,U_6\}$ itself.
By Proposition~\ref{ree5}, $\omega$ acts on the set $X$. Since $\omega$
interchanges the edges $\infty$ and $0$ and $U$
acts transitively on $X\backslash\{\infty\}$, we conclude that
(ii) and (iii) hold. Since $U_+$ is normal in $D_\infty$
(by \cite[4.7 and 5.3]{TW}) and $G$ is contained in
the centralizer of $\rho$, also (iv) holds.

For each $x\in X\backslash\{\infty\}$, there exists a unique
apartment $\Sigma_x$ of $\Gamma$ containing the edges $x$ and $\infty$.
For each vertex $u$ adjacent to $\bullet$ or $\star$,
there exists (by Proposition~\ref{ree6}) an edge $x\in X\backslash\{\infty\}$ such
that $u\in\Sigma_x$. If an element $g\in G$ acts trivially on 
$X$, then it acts trivially on all these apartments
and hence is itself trivial by \cite[3.7]{TW}. Thus (v) holds.
\end{proof}

\begin{proposition}\label{ree22}
Let $H$ be as in \eqref{ree1b}, let $D^\dagger$ be as in Notation~\ref{ree16},
let $D^\circ$ denote the centralizer of $\rho$ in $D^\dagger$, 
and let $T$ denote the two-point stabilizer 
$D^\circ_{\infty,0}$. Then there
is a canonical isomorphism $\pi$ from $H$ to 
$T$ that is compatible with the map $\varphi$ in Proposition~\ref{ree6}.
\end{proposition}

\begin{proof}
Let $g\in D^\dagger_{\infty,0}$. Thus $g$ acts trivially on the apartment $\Sigma$. 
By \cite[33.16]{TW}, there exist $a,u\in K^*$ such
that $x_1(s)^g=x_1(a^2 u^{-\theta}s)$ and $x_6(s)^g=x_6(a^{-\theta}u^2s)$ for all $s\in K$. 
By \eqref{ree23}, $g$ commutes with $\rho$ (and hence is contained in $T$)
if and only if $a^2u^{-\theta}=a^{-\theta}u^2$.
Since the maps $x\mapsto x^{2+\theta}$ and $x\mapsto x^{2-\theta}$ are
inverses of each other, we conclude that $a=u$ and the map
$g\mapsto a^{2+\theta}$ is an isomorphism from $T$ to $K^*$.
Now let $t=a^{2+\theta}$, so 
$x_1(s)^g=x_1(ts)$ and $x_6(s)^g=x_6(ts)$ for all $s\in K$. 
By the commutator relations~\eqref{ree10}, it follows that
$x_2(s)^g=x_2(t^{\theta+1}s)$, $x_3(t)^g=x_3(t^{\theta+2}s)$,
$x_4(s)^g=x_4(t^{\theta+2}s)$ and $x_5(s)^g=x_5(t^{\theta+1}s)$.
By Proposition~\ref{ree6}, therefore,
$(a,b,c)^g=h_t(a,b,c)$, where $h_t$ is as in \eqref{ree1b}. 
\end{proof}

\noindent
From now on we identify $H$ with the two point stabilizer $T$ via the map $\pi$ in Proposition~\ref{ree22}.

\smallskip
\section{The formula \eqref{ree1d}}\label{ree15}

In this section we show that 
the norm $N$ defined in \eqref{ree1c} is anisotropic
and that the automorphism $\omega$ satisfies \eqref{ree1d}.
We will do this by computing explicitly the action of $\omega$ on $X$
using Table~\ref{ta:m16}.

For each $g=(a,b,c)\in U$, we have
\begin{equation}\label{ree51a}
U_1 g = U_1\, x_2(b)\, x_3(c - ab + a^{\theta+2})\, x_4(c + ab)\, x_5(b - a^{\theta+1})\, x_6(a)
\end{equation}
by~Proposition~\ref{ree6} and 
\begin{equation}\label{ree51e}
U_1g\cap U=\{g\}
\end{equation}
by Proposition~\ref{ree7}(i).

\begin{lemma}\label{ree51b}
Suppose that $U_1x_2(\ddot{v})x_3(\ddot{u})x_4(\ddot{r})s_5(\ddot{s})x_6(\ddot{t})$
is one of the two vertices contained in an edge $\{U_1g,U_6g\}$ in the set
$X\backslash\{\infty\}$ for some $g\in U$. Then 
$g=(\ddot{t},\ddot{v},\ddot{r}-\ddot{v}\ddot{t})$.
\end{lemma}

\begin{proof}
This holds by \eqref{ree51a} and \eqref{ree51e}.
\end{proof}

We now fix $g=(a,b,c)\in U^*$ and let $u$, $v$ and $w=N(a,b,c)$ be as in Theorem~\ref{ree1}(iii).
Observe that the following curious identity holds:
\begin{equation}\label{eq:id3} 
w=av + bu + c^2.
\end{equation}
Let $m=m_1m_6$ (so $\omega=m^3$),
let $\alpha$ denote the vertex $U_1g$, let $\beta=\alpha^m$ and let $\gamma=\beta^m$.
Our goal is to show that $w\ne0$ and that
\begin{equation}\label{ree51d}
(a,b,c)^\omega=(-v/w,-u/w,-c/w)\, .
\end{equation}
 
\begin{lemma}\label{ree51c}
Suppose that $w\ne0$ and that
\[ \alpha^ \omega =
        U_1 \, x_2(\ddot{v}) \, x_3(\ddot{u}) \, x_4(\ddot{r}) \,
	x_5(\ddot{s}) \, x_6(\ddot{t}\,)\, . \]
Then \eqref{ree51d} holds if and only if
\begin{align}
\ddot{t} &= -v/w\, ;\label{eq:show1}\\
\ddot{v} &= -u/w\, ;\text{ and}\label{eq:show2}\\
\ddot{r} &= -c/w + (-v/w)(-u/w)\, .\label{eq:show3}
\end{align}
\end{lemma}

\begin{proof}
Since $\omega$ maps $X\backslash\{\infty,0\}$ to itself,
we have $\alpha^\omega=U_1e$ for some $e\in U^*$.
The claim holds, therefore, by Lemma~\ref{ree51b}.
\end{proof}

To begin, we assume that 
\begin{equation}\label{ree50a}
b\ne0\,, 
\end{equation}
so by Table~\ref{ta:m16} applied to \eqref{ree51a}, we have
\[ \beta=\alpha^m = U_1 \, x_2(\hat{v}) \, x_3(\hat{u}) \, x_4(\hat{r}) \,
	x_5(\hat{s}) \, x_6(\hat{t}\,) \,, \]
where
\begin{align}
\hat{v} &= b^{-1} c^\theta - a^\theta b^{\theta-1} + a^{2\theta+3} b^{-1} - c - ab \,; \notag\\
\hat{u} &= b - a^{\theta+1} + b^{-\theta} c^2 + a^2 b^{-\theta+2} + a^{2\theta+4} b^{-\theta} \notag\\
	& \hspace*{9ex} + a b^{-\theta+1} c - a^{\theta+2} b^{-\theta} c + 
          a^{\theta+3} b^{-\theta+1} \,; \notag\\[.3ex]
\hat{r} &= b^{-2} c^\theta - a^\theta b^{\theta-2} + a^{2\theta+3} b^{-2} + b^{-1}c - a \,; \notag\\
\hat{s} &= - b^{-\theta} c + a b^{-\theta+1} - a^{\theta+2} b^{-\theta} \,;\text{ and} \notag\\
\hat{t} &= b^{-1} \,.\label{ree60}
\end{align}
It is straightforward to check that the following identities hold:
\begin{align}
	w &= b \hat{u}^\theta - \hat{v} (\hat{v} - c) \,; \label{eq:id1}\\
	b \hat{r} &= \hat{v} - c \,; \label{eq:id2} \\
		b \hat{s}^\theta &= -a - b^{-1}(\hat{v}+c) \,;\text{ and} \label{eq:id4} \\
	\hat{v} &= -b^{-1} v \,. \label{eq:hatv}
\end{align}

Next we assume that
\begin{equation}\label{ree50b}
v\ne0\,.
\end{equation}
Thus also $\hat v\ne0$ (by~\eqref{eq:hatv}), so by a second application of 
Table~\ref{ta:m16}, we have
\[ \gamma=\beta^m = U_1 \, x_2(\tilde{v}) \, x_3(\tilde{u}) \, x_4(\tilde{r}) \,
	x_5(\tilde{s}) \, x_6(\tilde{t}) \,, \]
where
\begin{align}
	\tilde{v} &= \hat{u}^\theta \hat{v}^{-1} - \hat{r} \,; \label{ree60z}\\
	\tilde{u} &= \hat{s} + \hat{u}^2 \hat{v}^{-\theta} \,; \label{ree60x}\\
	\tilde{r} &= \hat{u}^\theta \hat{v}^{-2} + \hat{v}^{-1}\hat{r} + \hat{t} \,;\label{ree60y}\\
	\tilde{s} &= -\hat{u} \hat{v}^{-\theta} \,;\text{ and}\notag \\
	\tilde{t} &= \hat{v}^{-1}.\label{ree60b}
\end{align}
Note that
\begin{alignat}{2}
\tilde{v} &= \hat{u}^\theta \hat{v}^{-1} - \hat{r}& &\quad\text{by \eqref{ree60z}}\notag\\
&= \hat{v}^{-1} b^{-1} \cdot (b \hat{u}^\theta - b \hat{r} \hat{v})\notag\\
&=\hat v^{-1}b^{-1} \cdot \bigl(b\hat u^\theta-\hat v(\hat v-c)\bigr)& &\quad\text{by \eqref{eq:id2}}\notag\\
&= \hat{v}^{-1} b^{-1} \cdot w& &\quad\text{by \eqref{eq:id1}}\notag\\
&=-w/v& &\quad\text{by \eqref{eq:hatv}.}\notag
\end{alignat}
In particular, we have
\begin{equation}\label{eq:tildev2}
\tilde v=b^{-1}\hat v^{-1}w
\end{equation}
as well as 
\begin{equation}\label{eq:tildev3}
\tilde v=-wv^{-1}
\end{equation}
and $\hat u^\theta\hat v^{-1}=\hat r+b^{-1}\hat v^{-1}w$, so
\begin{equation}\label{eq:ida}
	b^2\hat{u}^{2\theta} \hat{v}^{-2} 
       = b^2\hat{r}^2 - b\hat{r} \hat{v}^{-1} w + \hat{v}^{-2} w^2\,.
\end{equation}
Moreover,
\begin{alignat}{2}
	\tilde{r}
&= \hat{u}^\theta \hat{v}^{-2} + \hat{v}^{-1}\hat{r} + \hat{t}& &\quad\text{by \eqref{ree60y}} \notag \\
&= \hat{u}^\theta \hat{v}^{-2} + \hat{v}^{-1}\hat{r} + b^{-1} & &\quad
            \text{by \eqref{ree60}} \notag \\
&= (\tilde v-\hat r)\hat v^{-1} + b^{-1}& &\quad\text{by \eqref{ree60z},}\notag
\end{alignat}
hence
\begin{equation}\label{eq:tilder}
\tilde r= b^{-1} \hat{v}^{-2} w - \hat{r} \hat{v}^{-1} + b^{-1} 
\end{equation}
by \eqref{eq:tildev2}, and thus 
\begin{equation}\label{ree71}
b\tilde rw=\hat v^{-2}w^2-b\hat{r}\hat{v}^{-1}w+w\,.
\end{equation}
We record also that
\begin{equation}\label{ree70}
b^2\hat s^\theta\hat v=-ab\hat v-\hat v^2-\hat vc
\end{equation}
by~\eqref{eq:id4}.

The vertex $\alpha^\omega=\gamma^m$ lies on an edge contained in $X$.
Hence $\alpha^\omega\in W_5$ (where $W_5$ is as in Figure~\ref{fig1}).
It follows that $\tilde v\ne0$
since otherwise $\gamma^m\in W_7$ by Table~\ref{ta:m16}.
By~\eqref{eq:tildev3}, we conclude that
$$w\ne0$$
and by a final application of Table~\ref{ta:m16}, we have
\[ \alpha^ \omega = \gamma^m = 
        U_1 \, x_2(\ddot{v}) \, x_3(\ddot{u}) \, x_4(\ddot{r}) \,
	x_5(\ddot{s}) \, x_6(\ddot{t}\,)\, ,  \]
where
\begin{align*}
	\ddot{v} &= \tilde{u}^\theta \tilde{v}^{-1} - \tilde{r} \,; \\
	\ddot{u} &= \tilde{s} + \tilde{u}^2 \tilde{v}^{-\theta} \,; \\
	\ddot{r} &= \tilde{u}^\theta \tilde{v}^{-2} + \tilde{v}^{-1}\tilde{r} + \tilde{t} \,; \\
	\ddot{s} &= -\tilde{u} \tilde{v}^{-\theta} \,;\text{ and} \\
        \ddot{t} &= \tilde{v}^{-1}\, .
\end{align*}
We now observe that $\ddot{t}=\tilde v^{-1}=-v/w$ by \eqref{eq:tildev3},
so \eqref{eq:show1} holds. Furthermore,
\begin{alignat}{2}
-b\ddot{v}w&=-b(\tilde u^\theta\tilde v^{-1}-\tilde r)w\notag\\
&=-b^2\tilde u^\theta \hat v+b\tilde rw& &\quad\text{by \eqref{eq:tildev2}}\notag\\
&=-b^2(\hat s+\hat u^2\hat v^{-\theta})^\theta\hat v+b\hat rw& &\quad\text{by \eqref{ree60x}}\notag\\
&=-b^2\hat s^\theta\hat v-b^2\hat u^{2\theta}\hat v^{-2}+b\tilde rw\, .\notag
\end{alignat}
Applying \eqref{eq:ida}, \eqref{ree71} and \eqref{ree70} to the three terms
in this last expression, we find that
\begin{alignat}{2}
-b\ddot{v}w&=ab\hat v+\hat v^2+c\hat v-b^2\hat r^2+w\notag\\
&=ab\hat v+\hat v^2+c\hat v-(\hat v-c)^2+w& &\quad\text{by \eqref{eq:id2}}\notag\\
&=-av-c^2+w& &\quad\text{by \eqref{eq:hatv}}\notag\\
&=bu& &\quad\text{by \eqref{eq:id3}}.\notag
\end{alignat}
Thus \eqref{eq:show2} holds. Finally, we have
\begin{alignat}{2}
	w \ddot{r}
	&= w ( \tilde{u}^\theta \tilde{v}^{-2} + \tilde{v}^{-1}\tilde{r} + \tilde{t}\, ) \notag\\
	&= w \tilde v^{-1}(\ddot{v}-\tilde r) + w \tilde{t} \notag\\
	&=-\tilde v^{-1}(u+w\tilde r)+w\tilde t& &\quad\text{by \eqref{eq:show2}}\notag\\
	&= u v w^{-1} + v \tilde{r} + w \tilde{t}& &\quad\text{by \eqref{eq:tildev3}}\notag \\
	&= u v w^{-1} + v\tilde r+w\hat v^{-1}
           & &\quad\text{by \eqref{ree60b}}\notag\\
	&= u v w^{-1} + v(b^{-1}\hat v^{-2}w-\hat v^{-1}\hat r+b^{-1})+w\hat v^{-1}
           & &\quad\text{by \eqref{eq:tilder}}\notag\\
	&= u v w^{-1} + (-\hat v^{-1}w+b\hat r-\hat v)+w\hat v^{-1}
           & &\quad\text{by \eqref{eq:hatv}}\notag\\
	&= u v w^{-1} - c & &\quad\text{by \eqref{eq:id2}}\,,\notag
\end{alignat}
so also \eqref{eq:show3} holds. By Lemma~\ref{ree51c}, it follows that
\eqref{ree51d} holds.
We conclude that $w\ne0$ and that the identity \eqref{ree1d} holds
for all ``generic'' points in $U^*$, i.e.~for all $(a,b,c)$ in $U^*$ satisfying
\eqref{ree50a} and \eqref{ree50b}. It is now
a much easier calculation to show using Table~\ref{ta:m16} that $w\ne0$ and
that the identity \eqref{ree1d} 
holds also when $b=0$ or $v=0$; we leave the details to the reader.

\smallskip
\section{Properties (I)--(VI)}

By Proposition~\ref{ree5}, $\omega$ is a permutation of $X$ of order~2. 
To conclude our proof of Theorem~\ref{ree1}, it thus
remains only to show that (I)--(VI) hold.
By Proposition~\ref{ree7}(ii) and (v), (I) holds. 
For each $x\in X$, there exists $g\in G$ mapping $\infty$
to~$x$; let $U_x=U^g$. If $g_1,g_2$ are two elements of
$G$ mapping $\infty$ to the same element of $X$, then
$g_1g_2^{-1}\in G_\infty$ and thus $U^{g_1}=U^{g_2}$ (since
$U_+$ is normal in $D_\infty$).
By Proposition~\ref{ree7}(i), it follows that
$(X,(U_x)_{x\in X})$ is a Moufang set (as defined, for 
example, in \cite[2.1]{tom}). Let $\mu$ be as in \cite[3.1]{tom}.
Thus for each $a\in U^*$, $\mu(a)$ is the unique element of
$U_0aU_0=U^\omega aU^\omega$ that interchanges $\infty$ and $0$. (Note that 
this is not the same $\mu$ as in the definition of $m_1$ and $m_6$
at the beginning of Section~\ref{ree17} above.) By \cite[3.1(ii)]{tom}, 
we have 
\begin{equation}\label{ree21a}
G_\infty=U\cdot\langle \mu(a)\mu(b)\mid a,b\in U^*\rangle \,.
\end{equation}

\begin{proposition}\label{ree21}
The following hold:
\begin{thmenum}
\item $G_\infty=UH^\dagger$, where $H^\dagger$ is as defined in \eqref{ree1e}.
\item $\omega\in\langle U,U^\omega\rangle$.
\end{thmenum}
\end{proposition}

\begin{proof}
We have $\langle \mu(a)\mu(b)\mid a,b\in U^*\rangle=H^\dagger$ by \cite[6.12(ii)]{petra}, whose proof 
depends only on knowing that the norm $N$ is anisotropic. 
By \eqref{ree21a}, therefore, (i) holds.
At the conclusion of the proof
of \cite[6.12(ii)]{petra}, it is observed that 
$\omega=\mu(0,0,1)$. Hence (ii) holds.
\end{proof}

By Propositions \ref{ree7}(iv) and \ref{ree21}, (II) and (III) hold.
Since $H$ normalizes both $U$ and $U^\omega$,
it follows from (III) that (IV) holds.
Let 
\begin{equation}\label{ree14b}
t\cdot(a,b,c)=h_t(a,b,c)
\end{equation}
for each $(a,b,c)\in U$ and each $t\in K^*$. By \eqref{ree1d}, we have
\begin{equation}\label{ree19}
\omega(t\cdot(a,b,c))=t^{-1}\cdot\omega(a,b,c)
\end{equation}
for all $(a,b,c)\in U$ and all $t\in K^*$. Thus (V) holds.

Suppose, finally, that $|K|>3$. Let $K^\dagger$ be as in \eqref{ree1e}.
Thus, in particular, $(K^*)^2=N(0,0,K^*)\subset K^\dagger$. Since $|K|>3$,
it follows that we can choose $t\in K^\dagger$ such that
$t^{\theta+1}\ne1$. Thus $t\ne1$, so also $t^{\theta+2}\ne1$. We have
\begin{align*}
	[h_t,(a,0,0)] &= \bigl( (1-t)a, (t-1)t^\theta a^{\theta+1}, 0 \bigr) \,, \\
	[h_t,(0,b,0)] &= \bigl( 0, (1-t^{\theta+1})b, 0 \bigr) \text{ and} \\
	[h_t,(0,0,c)] &= \bigl( 0, 0, (1-t^{\theta+2})c \bigr) \,
\end{align*}
for all $a,b,c\in K$.
Hence $U\subset[G,G]$. By Proposition~\ref{ree7}(iii), $(G_\infty,\langle\omega\rangle)$
is a BN-pair (as defined in \cite[2.1]{simple}). The group $U$ is nilpotent.
By \cite[2.8]{simple} and Proposition~\ref{ree7}(iv) and (v), it follows that
$G$ is simple. Thus (VI) holds.

\smallskip
\section{A more elementary reason why the norm is anisotropic}\label{ree11}

In this section we give a short algebraic proof that the norm $N$ defined in \eqref{ree1c} is
anisotropic. Let 
\begin{equation}\label{ree12}
\Omega(a,b,c)=(-v,-uw^\theta,-cw^{\theta+1})
\end{equation}
for all $(a,b,c)\in U$, where, as in \eqref{ree1c} and \eqref{ree1d},
\begin{align*}
v&=a^\theta b^\theta-c^\theta+ab^2+bc-a^{2\theta+3},\\
u&=a^2b-ac+b^\theta-a^{\theta+3}
\end{align*}
and $w=N(a,b,c)=-ac^\theta+a^{\theta+1}b^\theta-a^{\theta+3}b-a^2b^2+b^{\theta+1}+c^2-a^{2\theta+4}$. 
We first note that 
\begin{equation}\label{ree14}
N(\Omega(a,b,c))=N(a,b,c)^{2\theta+3}
\end{equation}
for all $(a,b,c)\in U$. This can be checked simply by 
plugging the definitions of $v$, $u$ and $w$ into \eqref{ree12}.
(That this identity {\it ought} to hold follows from \cite[6.18]{petra}
and~\eqref{ree19}.) 
Note, too, that 
\begin{equation}\label{ree14a}
N(t\cdot(a,b,c))=t^{2\theta+4}N(a,b,c) 
\end{equation}
and 
\begin{equation}\label{ree14d}
N\big((a,b,c)^{-1}\big)=N(a,b,c)
\end{equation}
for all $(a,b,c)\in U$ and all $t\in K^*$, where
$t\cdot(a,b,c)$ is as in \eqref{ree14b}
and $(a,b,c)^{-1}$ is as in Theorem~\ref{ree1}(i).

Now fix $(a,b,c)\in U^*$ such that $w=0$. 

\begin{lemma}\label{ree80}
$v=0$.
\end{lemma}

\begin{proof}
By \eqref{ree12} and \eqref{ree14}, we have
\[ N(-v,0,0)=N(\Omega(a,b,c))=0 \,. \]
By \eqref{ree1c}, on the other hand, $N(-v,0,0) = -v^{2\theta+4}$. 
\end{proof}

\begin{lemma}\label{ree81}
$a\ne0$.
\end{lemma}

\begin{proof}
Suppose $a=0$. Since $(a,b,c)\ne0$ and $w=0$, we have $c\ne0$. By \eqref{ree14a}, the norm of 
$c^{\theta-2}\cdot(0,b,c)$ is zero. We can thus assume that $c=1$. It follows by
\eqref{ree1c} that $b\ne1$ but by Lemma~\ref{ree80} that $b=1$.
\end{proof}

By \eqref{ree14a} and Lemma~\ref{ree81}, we can assume from now on that 
$a=1$. Hence $v=0$ means that 
\begin{equation}\label{ree84}
b^\theta-c^\theta+b^2+bc-1=0
\end{equation}
and $w-v=0$ means that
\begin{equation}\label{ree85}
b^{\theta+1}+b^2-b-bc+c^2=0 \,.
\end{equation}
By \eqref{ree14d} and Lemma~\ref{ree80}, we also have $v(-1,-b+1,-c)=v\big((1,b,c)^{-1}\big)=0$
and thus 
\begin{equation}\label{ree86}
b^\theta+c^\theta-b^2-b-1+bc-c=0 \,.
\end{equation}
Adding \eqref{ree84} and \eqref{ree86}, we find that
\begin{equation}\label{ree87}
	b^\theta + b - 1 = - bc - c \,.
\end{equation}
Multiplying this last equation by $b$ and comparing with \eqref{ree85}, we obtain
\begin{equation}\label{ree88}
	c(c - b^2 + b) = 0 \,.
\end{equation}
Assume first that $c=0$.
Then by \eqref{ree84}, we have $b^\theta + b^2 - 1 = 0$ whereas by \eqref{ree87},
we have $b^\theta + b - 1 = 0$.
We find $b^2 = b$ and thus $b \in \{ 0,1 \}$, contradicting the equality $b^\theta + b - 1 = 0$.

Hence $c \neq 0$, and it follows from \eqref{ree88} that $c = b^2 - b$.
By \eqref{ree84}, we now obtain
\[ b^{2\theta} = b^3 - 1 - b^\theta \,; \]
from \eqref{ree87} on the other hand, we get
\[ b^3 - 1 = - b^\theta \,. \]
Combining the last two equations, we obtain $b^{2\theta} = b^\theta$,
but then $c^\theta = 0$ and hence $c=0$ after all.
With this contradiction, we conclude that the norm $N$ is anisotropic.

\smallskip
\section{The subgroup $H^\dagger$}\label{ree18}

If $K$ is finite, then $|K|$ is an odd power of 3, from which it follows that 
$K^*$ is generated by $(K^*)^2=N(0,0,K^*)$ and $-1=N(0,1,1)$, 
so $K^\dagger=K^*$ and $H^\dagger=H$. This is \cite[8.4]{ree}.
It is not necessarily true, however, that $H^\dagger=H$ if $K$ is infinite.
In this section we illustrate this with
an example. As Tits suggests in \cite[1.12]{octagons}, we only 
need to modify what he does there slightly.

Let $F$ be an odd degree extension of the field with three elements
and let K be the field of quotients of the polynomial ring
$F[s,t]$ in two variables $s$ and $t$. Since $|F|$ is an odd power of~3, there
exists a unique endomorphism 
$\theta$ of $K$ mapping $F$ to $F$, $t$ to $s$ and $s$ to $t^3$ 
whose square is the Frobenius endomorphism.
(In what follows, the reader may wish to think of $s$ as being formally equal to $t^{\sqrt{3}}$.)

\begin{proposition}\label{ree30}
The group $K^\dagger\cap F(t)$ is generated by $(F(t)^*)^2$ and 
all irreducible polynomials in $F[t]$ of even degree.
\end{proposition}

\begin{proof}
Since $F$ is finite, we have $F^*\subset K^\dagger$.
Let $f\in F[t]$ be an irreducible polynomial of even degree over $F$
and let $\alpha$ be a root of $f$ in some splitting field $L$. 
Then $L=F(\alpha)$ and $[L:F]={\rm deg}(f)=2d$ for some $d$. Thus
$L$ contains an element $\beta$ whose square is $-1$. Since
$[L:F(\beta)]=d$, there are non-zero polynomials $p,q\in F[t]$ of degree
at most $d$ such that $p+\beta q$ is the minimal polynomial of
$\alpha$ over $F(\beta)$. Thus $p+\beta q$ divides $f$. 
Hence also $p-\beta q$ divides $f$. Since the polynomial $p+\beta q$ is irreducible
over $F(\beta)$, it follows that it is relatively prime to the 
polynomial $p-\beta q$. Thus $f/e$ equals 
the product of these two polynomials for some $e\in F^*$. Hence
\[ f=e(p^2+q^2)=eN(0,p^{\theta-1},q)\in K^\dagger \,. \]

Now suppose that $g\in F[t]$ is the product of distinct irreducible polynomials
of odd degree. It will suffice to show that $g\not\in K^\dagger$. 
Let $F_1$ be the splitting field of $g$ over $F$ and let $K_1=F_1(s,t)$. 
The extension $F_1/F$ is of odd degree by the choice of $g$, so $\theta$ has
a unique extension to an endomorphism of $K_1$ (which we continue to call
$\theta$) whose square is the Frobenius map.
Let $c$ be an arbitrary root of $g$ in $F_1$ and let $d=c^\theta$. We define a valuation
$\nu$ on $K_1$ with values in ${\mathbb Z}[\sqrt{3}]$.
First we declare the degree of a monomial $e(s-d)^m(t-c)^n$ (for $e\in F_1^*$) to be 
$n+m\sqrt{3}$. If $p\in F_1[s,t]^*$, we write $p$ as a sum of monomials
in the variables $t-c$ and $s-d$ and define $\nu(p)$ to be the minimum
of the degrees of these monomials (minimum with respect to the natural ordering of 
${\mathbb Z}[\sqrt{3}]$ as a subset of ${\mathbb R}$). Finally we set
$\nu(p/q)=\nu(p)-\nu(q)$ for all $p,q\in F_1[s,t]$. Then $\nu$ is a
well defined valuation on $K_1$ such that $\nu(g)=1$ and $\nu(a^\theta)=\sqrt{3}\cdot\nu(a)$
for all $a$. 

Now let $w=N(a,b,c)$ for $a,b,c\in K_1$. 
By \cite[9.3]{petra} (whose proof depends only on the fact that the
norm is anisotropic), $\nu(w)$ is equal to the minimum of
$(2\sqrt{3}+4)\nu(a)$, $(\sqrt{3}+1)\nu(b)$ and $2\nu(c)$. 
Since $(\sqrt{3}+1)^2 = 2\sqrt{3}+4$ and $(\sqrt{3}+1)(\sqrt{3}-1)=2$, it follows that
$\nu(K_1^\dagger)=(\sqrt{3}+1){\mathbb Z}[\sqrt{3}]$.
Since $\nu(g)=1\not\in(\sqrt{3}+1){\mathbb Z}[\sqrt{3}]$, we conclude that
$g\not\in K_1^\dagger$. Hence $g\not\in K^\dagger$.
\end{proof}

\begin{corollary}\label{ree31}
$K^*/K^\dagger$ is infinite.
\end{corollary}

\begin{proof}
There are infinitely many pairwise non-proportional irreducible polynomials
of odd degree in $F[t]$. By Proposition~\ref{ree30}, these polynomials have pairwise
distinct images in $K^*/K^\dagger$.
\end{proof}

\end{document}